\documentclass[11pt]{article}
\newcommand{\newsection}[1]{\setcounter{equation}{0}\section{#1}}

\def\keywords{ \if@twocolumn
\section*{Keywords}
\else \small
\begin{center}
{ \bf Keywords\vspace{-.5em}\vspace{0pt}}
\end{center}
\center \fi}
\def\endkeywords{ \if@twocolumn\else\endcenter\fi}
\usepackage{epsfig}
\usepackage{graphicx}
\begin{document}
\title{\bf $N$ coupled non-local harmonic oscillators leading to
$2N$-th Order Initial Value Problem}
\author{\thanks{
Department of Mathematics, University of the Punjab, Lahore 54590,
 Pakistan.   ~~~~~~~~~~~~~~~~~~~~~~~~~~~~~~~~~~~~Email: toghazala2003@yahoo.com} Ghazala Akram \thanks{
School of Mathematical Sciences, Queen Mary, University of London,
Mile End Road, London E1 4NS, UK. ~~Email:
 g.akram@qmul.ac.uk} ,
 ~Christian Beck \thanks{School of Mathematical Sciences, Queen Mary, University of London,
Mile End Road, London E1 4NS, UK.~~~Email:
 c.beck@qmul.ac.uk}}
\date{}
\maketitle
\begin{center}
\begin{minipage}{5.0in}
\begin{abstract}
\noindent We consider a set of interwoven harmonic oscillators where
the acceleration of a given oscillator is determined by the position
of its nearest neighbor. We show that this problem of $N$ non-local
oscillators with periodic boundary conditions leads to a $2N$-th
order initial value problem. We discuss the numerical solution of
this using a non-polynomial spline method. A very precise numerical
method that minimizes the error can be developed, which we test for
a few examples of driving forces.
\end{abstract}
\end{minipage}
\end{center}
{\bf Keywords:} {\small Non-local coupled harmonic oscillators, non-polynomial spline; consistency
relations; end conditions; $N$th-order Initial Value Problem.}
\newsection{Introduction}

Coupled harmonic oscillators are a standard paradigm in
many engineering, physical, chemical or biological systems. The basic ingredient
of the harmonic oscillator dynamics is the fact that the acceleration
of each oscillator
is proportional to its position, with a proportionality constant
given by the negative of the frequency squared. In addition, for
non-equilibrium situations there are
time-dependent driving forces on each oscillator.

In this paper we are interested in a fundamental modification of
this approach. We assume that the acceleration of each oscillator
is proportional to the position of its neighboring oscillator.
By this, of course, a strongly coupled structure is introduced,
and a kind of non-local dynamics,
 since the acceleration of a given oscillator is determined by a position variable elsewhere.
 Such a dynamics is  motivated for very 
strongly coupled individual systems.
For example in \cite{os1,os2,os3} 
strongly coupled oscillator systems are studied,
which degenerate to our type of dynamics in the infinite coupling limit.
In this limit individual local oscillator behavior is influenced in a 
hierarchical way by nearest neighbors.
We assume a driven non-equilibrium situation where each oscillator is also driven by
individual time-dependent driving forces.

We will show that the above dynamics, for $N$ oscillators with
periodic boundary conditions, leads to a high-order initial value
problem, indeed of order $2N$. Thus rather large derivatives will
become relevant if there are just a few oscillators coupled in this
way. These types of initial value problems require effective
numerical methods, which depend on the number $N$. Useful in this
context are non-polynomial spline methods. We will describe the
optimum way to solve this system numerically for a given set of
driving forces. In fact, we will show that there is an optimum
method with minimum error, which is in particular useful if highest
precision numerical results are required. As main examples, we will
deal with the cases $N=2$ and $N=3$.

\section{Nonlocal coupling of driven harmonic oscillators}

To illustrate the idea we start with two uncoupled driven harmonic oscillators:

\begin{eqnarray}
\ddot{y_1}+\omega_1^2  y_1&=&g_1 (t) \\
\ddot{y_2}+\omega_2^2  y_2&=&g_2(t)
\end{eqnarray}
$\omega_i$ is the frequency of oscillator $i$,
and $g_i(t)$ is a local driving force on oscillator $i$.
As mentioned in the introduction, the basic idea is to consider a modified
dynamics where the position $y_1$ of the first oscillator determines the
acceleration force of the second one, and vice versa. That is to say,
instead of the above trivial (uncoupled) dynamics we consider the following coupled dynamics:
\begin{eqnarray}
\ddot{y_1}+\omega_1^2 y_2&=&g_1(t) \label{111}\\
\ddot{y_2}+\omega_2^2 y_1&=&g_2(t) \label{222}
\end{eqnarray}
Differentiating eq.~(\ref{222}) twice we get
\begin{equation}
y_2^{(4)}+\omega_2^2 y_1^{(2)}=g_2^{(2)}(t)
\end{equation}
where $y^{(k)}$ denotes the $k$-th derivative with respect to time $t$.
Eliminating in this equation $y_1^{(2)}$ using eq.~(\ref{111}) we
get
\begin{equation}
y_2^{(4)}+\omega_2^2(g_1(t)-\omega_1^2y_2)=g_2^{(2)} (t)
\end{equation}
This is equivalent to a 4-th order equation of the form
\begin{equation}
y_2^{(4)}+f(t)y_2(t)=g(t)
\end{equation}
where
\begin{equation}
f(t)=-\omega_1^2\omega_2^2
\end{equation}
and
\begin{equation}
g(t)=g_2^{(2)}(t)-\omega_2^2 g_1(t).
\end{equation}
Together with the initial position and velocity of the two oscillators,
this leads to an initial value problem of 4th-order, for which we will describe
the optimum numerical method in section 3.

But let us here first extend the
problem, by considering $N=3$ non-local oscillators, and later an arbitrary
number $N$. The non-locally coupled dynamics for 3 oscillators reads
\begin{eqnarray}
\ddot{y_1}+\omega_1^2 y_2&=&g_1(t) \label{1111}\\
\ddot{y_2}+\omega_2^2 y_3&=&g_2(t) \label{2222}\\
\ddot{y_3}+\omega_3^2 y_1&=&g_3(t) \label{3333}
\end{eqnarray}
Differentiating eq.~(\ref{3333}) twice we get
\begin{equation}
y_3^{(4)}+\omega_3^2 y_1^{(2)}=g_3^{(2)}(t).
\end{equation}
Eliminating in this equation $y_1^{(2)}$ using eq.~(\ref{1111}) we
get
\begin{equation}
y_3^{(4)}+\omega_3^2(g_1(t)-\omega_1^2y_2)=g_3^{(2)} (t)
\end{equation}
or
\begin{equation}
y_3^{(4)}-\omega_1^2\omega_3^2y_2=g_3^{(2)}(t)-\omega_3^2g_1(t).
\end{equation}
Differentiating twice this leads to
\begin{equation}
y_3^{(6)}-\omega_1^2\omega_3^2y_2^{(2)}=g_3^{(4)}(t)-\omega_3^2g_1^{(2)}(t)
\end{equation}
and using (\ref{2222}) to eliminate $y_2^{(2)}$ we get
\begin{equation}
y_3^{(6)}+\omega_1^2\omega_2^2\omega_3^2y_3=g_3^{(4)}(t)-\omega_3^2g_1^{(2)}(t)+\omega_3^2\omega_1^2g_2(t).
\end{equation}
Apparently, this can be written in the form
\begin{equation}
y_3^{(6)}+f(t)y_3(t)=g(t)
\end{equation}
where
\begin{equation}
f(t)=+\omega_1^2\omega_2^2\omega_3^2
\end{equation}
and
\begin{equation}
g(t)=g_3^{(4)}(t)-\omega_3^2 g_1^{(2)}(t)+\omega_3^2\omega_1^2g_2(t).
\end{equation}

It is obvious how to generalize this problem to $N$ oscillators.
In this case one obtains
a $2N$-th order initial value problem of the form
\begin{equation}
y_N^{(2N)}+f(t) y_N=g(t)
\end{equation}
where
\begin{equation}
f(t)=(-1)^{N+1} \prod_{i=1}^N \omega_i^2
\end{equation}
and $g(t)$ is a sum of various derivatives of local driving forces $g_i(t)$ weighted with
frequencies.

The optimum numerical way to solve these high-order
initial value problems depends on $N$ in a nontrivial way. In the following, we allow
for general time-dependent functions $f(t)$ and deal in
detail with the cases $N=2$ and $N=3$.

\section{N=2: Optimized numerical solution of the 4th order initial value problem}

\subsection{The initial value problem}
The case of two oscillators ($N=2$) leads to the following fourth order initial value problem
\begin{equation}\label{e:1.1}
\left. \begin{array}{lcl}
    y^{(4)}(t)+f(t)y(t)&=&g(t),  \ \ \ \ t \ \in \  [a, \ b], \\
    y(a)  \ \ \ = u_{0},\ & &  y^{(1)}(a) \ \ \  =  u_{1}, \\
    y^{(2)}(a) = u_{2},\ & & y^{(3)}(a)  =   u_{3},
    \end{array} \right \}
\end{equation}
where the $ u_{i} (i= 0,1,2,3)$ are finite real constants while the
functions $f(t)$ and $g(t)$ are continuous on $[a, b]$. To simplify the notation,
we have written $y(t)$ instead of $y_2(t)$.

Of course, from an engineering point we are interested in four initial values
given by initial position and velocity of oscillator 1 and 2. These are related
to the constants $u_i$ by
\begin{eqnarray}
y_2(a)&=&u_0 \\
\dot{y_2} (a)&=& u_1 \\
y_1(a)&=&\frac{g_2(a)-u_2}{\omega_2^2} \\
\dot{y_1} (a)&=&\frac{\dot{g_2}(a)-u_3}{\omega_2^2}
\end{eqnarray}
We have reduced the 4-dimensional problem of the space-space structure of the two coupled nonlocal
oscillators to a 1-dimensional initial value problem of 4-th order, for which we can apply very
precise numerical methods, as described in the following subsection. While for our nonlocal
oscillator example $f(t)=const=-\omega_1^2\omega_2^2$, the numerical method developed in the following is
applicable general time-dependent functions $f(t)$, as long as they are continuous.

  \subsection{Nonpolynomial Spline Method}
  \indent To develop the spline approximation to the problem $(\ref{e:1.1})$,
  the interval $[a, \ b]$ is divided
  into $n$ equal subintervals, using the grid points $t_{i}  =  a +  ih $ ;
   $i=0,1,\ldots ,
  n$, where $h \ = \ (b-a)/n.$ \\
    Consider the following restriction $S_{i}$ of the solution to each subinterval $[t_{i}, \
  t_{i+1}], \ i=0,1,\ldots , n-1$,
  \begin{eqnarray}
  S_{i}(t)& = &a_{i}\cos \omega (t-t_{i})+b_{i}\sin \omega (t-t_{i})+c_{i}(t-t_{i})^{3}
  +d_{i}(t-t_{i})^{2}+e_{i}(t-t_{i})+ p_{i} \nonumber \\
  & & ~~~~~~
  \end{eqnarray}
     Let
  \begin{equation}
 \left. \begin{array}{rcl}
  y_{i} \ = \ S_{i}(t_{i}) & & \ \ \ \ M_{i} \ = \ S_{i}^{(2)}(t_{i}),
    \\
  N_{i} \ = \ S_{i}^{(4)}(t_{i}),  & &
  \end{array} \right \} \ i=0,1,\ldots,n.
  \end{equation}
  Following~\cite{Sg08} and postulating that at the end points of the intervals the 1st and 3rd  derivatives
  are continuous, one derives
  the following consistency relation between the values of splines and their fourth order derivatives at
  border points:
\begin{eqnarray} \label{e:2.8}
& \left(\alpha h^{4}N_{i-4}+\beta h^{4}N_{i-3}+\gamma
h^{4}N_{i-2}+ \beta h^{4}N_{i-1}+\alpha h^{4}N_{i} \right) \nonumber \\
= & \left[y_{i-4} -4 y_{i-3} +6 y_{i-2} -  4 y_{i-1}+ y_{i}
   \right]; \nonumber \\
  & i=4,5,\ldots,n,
\end{eqnarray}
where
\begin{eqnarray}
\alpha  =  \left(\frac{1}{6 \theta \sin \theta}- \frac{1}{
\theta^{3} \sin \theta}+\frac{1}{ \theta^{4}}\right), ~~~~~ \beta =
 \left(\frac{2(1+\cos \theta)}{\theta^{3}\sin \theta}-
\frac{( \cos \theta -2)}{3 \theta}-\frac{4}{ \theta^{4}} \right)
\nonumber
\end{eqnarray}
and
\begin{eqnarray}
 \gamma  =  \left(\frac{-2(1+2 \cos \theta)}{\theta^{3}\sin \theta}+
\frac{(1-4 \cos \theta )}{ 3 \theta \sin \theta}+\frac{6}{
\theta^{4}} \right). & & \nonumber
\end{eqnarray}
Here $\theta=\omega h$ is an arbitrary parameter. The relation
$(\ref{e:2.8})$ forms a system of $n-3$ linear
equations in the $n$ unknowns \\
$(y_{i}, \ i=1,2,...,n) $, while $N_{i}$ is taken from IVP
$(\ref{e:1.1})$ to be equal to $-f_{i}y_{i} + g_{i}$, \\ $i=0,1,\ldots,n.$ \\
 Following~\cite{Sg074}, three equations (end conditions) are
determined to find the complete solution of $y_{i}s$ appearing in
eq. $(\ref{e:2.8})$, as given below:
\begin{eqnarray}\label{e:3.1}
  N_{0}  +  N_{4} & = & \frac{1}{h^{4}} \ \left[  \ \frac{-220}{9} y_{0}
   +40 y_{1} -20 y_{2} +  \frac{40}{9}y_{3} -  \frac{40}{3}hy_{0}^{(1)}
     - \frac{4}{3}h^{4}y_{0}^{(4)}
     \right ], \\ \label{e:3.2}
     N_{1}  +  N_{5} & = & \frac{1}{h^{4}} \ \left[  \ \frac{18336}{575} y_{1}
   -  \frac{22992}{575} y_{2} +  \frac{4656}{575} y_{3} + \frac{2736}{115}hy_{0}^{(1)}
          \right. \nonumber \\
    & & \ \ \ \ \ \ \ \ \ \left.  +  \frac{15864}{575} \ h^{2}y_{0}^{(2)}
    +  \frac{6648}{575} \ h^{3}y_{0}^{(3)}  \right ]
   \end{eqnarray}
and
\begin{eqnarray}\label{e:3.3}
  N_{2}  +   N_{6}  & = & \frac{1}{h^{4}} \ \left[ \ \frac{8157}{865}y_{2}
   -  \frac{11424}{865}y_{3} +  \frac{3267}{865} y_{4}
        +  \frac{978}{173}hy_{0}^{(1)} +   \frac{8958}{865} \ h^{2}y_{0}^{(2)} \right. \nonumber \\
   & & \left.
    + \frac{5684}{865}h^{3}y_{0}^{(3)}     \right]\ .
    \end{eqnarray}
    The local truncation errors associated with the linear equations
    $(\ref{e:3.1})-(\ref{e:3.3})$ and $(\ref{e:2.8})$ are calculated as
\begin{equation}\label{eee:3.12}
\tilde{t}_{i} = \left\{ \begin{array}{ll}
   - \frac{47}{9} h^{6}y^{(6)}(t_{1})+O(h^{7}), &  i=1,  \\
   - \frac{71686}{8625} h^{6}y^{(6)}(t_{2})+O(h^{7}), &  i=2,  \\
    - \frac{143342}{12975} h^{6}y^{(6)}(t_{3})+O(h^{7}), &  i=3,  \\
   \frac{1}{6}(-1+24 \alpha +6 \beta) h^{6}y^{(6)}(t_{i})+O(h^{7}), &  i=4,5,\ldots,n  \\
        \end{array} \right.
\end{equation}
and
\begin{equation}\label{eee:3.13}
\|\tilde{T}\|=ch^{6}R_{1}=O(h^{6}), \ \ \ \ \ \ \ \ R_{1}= \max_{t
\in [a, \ b]} |y^{(6)}(t)|,
\end{equation}\
 where $c$ is a constant which
depends only on the values of $\alpha$ and
$\beta$ and is independent of $h$. \\

Let us mention that the solution obtained using the system of linear
equations $(\ref{e:3.1})-(\ref{e:3.3})$ and $(\ref{e:2.8})$ is
second order convergent. But if $\alpha$, $\beta$ and $\gamma$ are
taken such that $\alpha=-
\frac{1}{720},~~\beta=\frac{31}{180},~~\gamma=\frac{79}{120}$ then
the order of the truncation error in eq. $(\ref{e:2.8})$ is
$O(h^{10})$ and the order of convergence can then be improved up to
sixth order, using this method of improved order of end conditions.\\ \\
The improved end conditions with truncation error of order $O(h^{10})$
are
\begin{eqnarray}\label{e:3.4}
  & N_{0} +\frac{843268}{2081}N_{1}+\frac{330342}{2081}N_{2}-\frac{16892}{2081}N_{3} +  N_{4} \nonumber \\
   =
  & \frac{1}{h^{4}} \ \left[  \ \frac{-68397280}{18729} y_{0}
   +\frac{13366080}{2081} y_{1}  -\frac{7408800}{2081} y_{2} +  \frac{14781760}{18729}y_{3}
   -  \frac{10427200}{6243}hy_{0}^{(1)}+\frac{743680}{2081}h^{2}y_{0}^{(2)}
      \right. \nonumber \\
    &  \ \ \ \ \ \ \ \ \ \left.+ \frac{259840}{2081}h^{3}y_{0}^{(3)}
     \right ],  \\ \label{e:3.5}
    & N_{1} -\frac{156090207332}{158360705}N_{2}-\frac{40456201386}{158360705}N_{3}-\frac{600708692}{158360705}N_{4}
      +  N_{5} \nonumber \\
       = & \frac{1}{h^{4}} \ \left[  \ \frac{180155114496}{31672141} y_{1}
   -  \frac{340726283352}{31672141} y_{2} + \frac{210168798336}{31672141} y_{3}-  \frac{49597629480}{31672141} y_{4} + \frac{69181575120}{31672141}hy_{0}^{(1)}
         \right. \nonumber \\
     & \ \ \ \ \ \ \ \ \ \left.
    +  \frac{42396452784}{31672141} \ h^{2}y_{0}^{(2)} +  \frac{7557647328}{31672141} \ h^{3}y_{0}^{(3)}  \right ]
   \end{eqnarray}
and
\begin{eqnarray}\label{e:3.6}
 &  N_{2} - \frac{85514900495708}{1252977040745}N_{3}+\frac{3759590586966}{1252977040745}N_{4}
  -\frac{7418340285788}{1252977040745}N_{5} +   N_{6} \nonumber \\
   =
  & \frac{1}{h^{4}} \ \left[ \ \frac{43463161469952}{250595408149}y_{2}
   -  \frac{94491207986112}{250595408149}y_{3}
   +  \frac{68699611790208}{250595408149} y_{4} -  \frac{17671565274048}{250595408149} y_{5} \right. \nonumber \\
    & \left.
        +  \frac{10106680227840}{250595408149}hy_{0}^{(1)}
        +   \frac{9581784601536}{250595408149} \ h^{2}y_{0}^{(2)}
    + \frac{2621304758016}{250595408149}h^{3}y_{0}^{(3)}     \right]\ .
    \end{eqnarray}
    The truncation errors of the corresponding equations are
\begin{equation}\label{e:3.12}
\tilde{t}_{i} = \left\{ \begin{array}{ll}
   - 0.3034 h^{10}y^{(10)}(t_{1})+O(h^{11}), &  i=1,  \\
   - 1.4034 h^{10}y^{(10)}(t_{2})+O(h^{11}), &  i=2,  \\
    - 1.0163 h^{10}y^{(10)}(t_{3})+O(h^{11}), &  i=3,  \\
   \frac{1}{30240}(-17+5376 \alpha +84 \beta) h^{10}y^{(10)}(t_{i})+O(h^{11}), &  i=4,5,\ldots,n  \\
        \end{array} \right.
\end{equation}
and
\begin{equation}\label{e:3.13}
\|\tilde{T}\|=ch^{10}R_{2}=O(h^{10}), \ \ \ \ \ \ \ \ R_{2}= \max_{t
\in [a, \ b]} |y^{(10)}(t)|,
\end{equation}\
 where $c$ is a constant which
depends only on the values of $\alpha$ and
$\beta$ and is independent of $h$. \\
\\
 To illustrate the powerfulness of the method, two
analytically solvable examples are discussed in the following:
 \subsection{Examples}
   \ \ \\
{\large {\bf Example 1 }} \\
  Consider the following IVP
 \begin{equation}
\left.
\begin{array}{rl}\label{e:4.1}
    y^{(4)}(t) \ - \  y(t) \ & =  \  4 cos(t), \ \ \ \ \ t \in [-1, \ 1],  \\
 y(-1) \ & =  -2 sin(1), \\
   y^{(1)}(-1) \ & =  \ 2cos(1)+sin(1), \\
   y^{(2)}(-1) \ & =  -2cos(1)+2sin(1)\\
   y^{(3)}(-1) \ & = -2cos(1)-3sin(1).
\end{array}
\right\}
\end{equation}
The analytic solution of the above problem is
$$ y(t) \ = \ (1-t) \  sin(t) \ .$$
The observed maximum errors (in absolute values) associated with $
y_{i}$, for the problem $(\ref{e:4.1})$, corresponding to the
different values of $\alpha, \ \beta$ and $\gamma$, are summarized
in Table 1. It is confirmed from Table 1 that if $h$ is reduced by
factor $1/2$, then $\| E \|$ is reduced by a factor $1/4$, which
indicates that the method gives second-order results.
\par \noindent
\begin{table}[htp]
\caption{ Maximum absolute errors for problem $(\ref{e:4.1})$ in
 $ y_{i}.$ }
 \begin{center}
\begin{tabular}{|c|c|c|c|} \hline
               &                              &
               &                                           \\
$n$              &  $\alpha=0, \ \beta=0 $
                &
                  $\alpha=1/2, \ \beta=1/2 $  &
                  $\alpha=1/6, \ \beta=1/6 $
                    \\
               &               &
               &
                                                       \\
               &  $\gamma=1 $             &
              $\gamma=-1 $ &
              $\gamma=1/3 $
                                                      \\
               &                &
               &
                                                            \\ \hline
               &                              &
               &                                            \\
$6 $        & $ 6.74 \times {10}^{-1} $
               & $ 3.6 \times {10}^{0} $
               & $ 1.73 \times {10}^{0} $
                                                                      \\ \hline
               &                              &
               &                                                                                 \\
$ 12$   & $ 5.77 \times {10}^{-2} $
               & $ 7.3 \times {10}^{-1} $
               & $ 2.22 \times {10}^{-1} $
                                                                                     \\ \hline
               &                              &
               &                                                                                 \\
$24 $  & $ 3.3 \times {10}^{-3} $
               & $ 4.5 \times {10}^{-2} $
               & $ 1.3 \times {10}^{-2} $
                                                                  \\ \hline
               &                              &
               &                                                                                \\
$48$ & $ 1.48 \times {10}^{-4} $
               & $ 2.1 \times {10}^{-3} $
               & $ 5.93 \times {10}^{-4} $
                                               \\ \hline
 \end{tabular}
\end{center}
\end{table}
\newpage
The observed maximum errors (absolute values) associated with $
y_{i}$ for the problem $(\ref{e:4.1})$, corresponding to the
use of improved end conditions, are summarized in Table 2.
\par \noindent
\begin{table}[htp]
\caption{ Maximum absolute errors for problem $(\ref{e:4.1})$ in
 $ y_{i}.$ }
 \begin{center}
\begin{tabular}{|c|c|} \hline
               &                                          \\
$n$        & $|y(t_{i})-y_{i} |$
                                                            \\ \hline
               &                                            \\
$6$       & $ 1.7 \times {10}^{-3} $
                                                                                     \\ \hline
               &                                                                               \\
$ 12$   & $ 1.17 \times {10}^{-5} $
                                                                                                    \\ \hline
               &                                                                                \\
$24 $  & $ 7.19 \times {10}^{-8}$
                                                                  \\ \hline
                                                                   &                                                                                \\
$48 $  & $ 7.72 \times {10}^{-11}$
                                                                  \\ \hline
 \end{tabular}
\end{center}
\end{table}
A significant improvement of precision is obtained.

\par \noindent
{\large {\bf Example 2 }} \\
 \ \ \\
 Consider the following IVP
 \begin{equation}
\left.
\begin{array}{ll}\label{e:4.2}
    y^{(4)}(t) \ + t  y(t) \ =  \ - e^{t}(8+7t+t^{3}), \ \ \ \ \ 0\leq t \leq 1 \\
 y(0) \  \ \ \  = 0,  \ \ \ \ \ \ \ \ \ \ \ \ \ \ \ \  y^{(1)}(0) \ = 1, \\
   y^{(2)}(0) \ =  0, \ \ \ \ \ \ \ \ \ \ \ \ \ y^{(3)}(0) \ =  -3.
\end{array}
\right\}
\end{equation}
The corresponding analytic solution is now
$$ y(t) \ = \ t(1-t) \  e^{t} \ .$$
The observed maximum errors for different values of $\alpha, \
\beta$ and $\gamma$ are summarized in Table 3. It is confirmed from
Table 3 that if $h$ is reduced by factor $1/2$, then $\| E \|$ is
reduced by a factor $1/4$, which indicates that the method gives
second-order results.
\par \noindent
\begin{table}[htp]
\caption{ Maximum absolute errors for problem $(\ref{e:4.2})$ in
 $ y_{i}.$ }
 \begin{center}
\begin{tabular}{|c|c|c|c|} \hline
               &                              &
               &                                           \\
$n$              &  $\alpha=0, \ \beta=0 $
                &
                  $\alpha=1/2, \ \beta=1/2 $  &
                  $\alpha=1/6, \ \beta=1/6 $
                    \\
               &               &
               &
                                                       \\
               &  $\gamma=1 $             &
              $\gamma=-1 $ &
              $\gamma=1/3 $
                                                      \\
               &                &
               &
                                                            \\ \hline
               &                              &
               &                                            \\
$6 $        & $ 1.14 \times {10}^{-1} $
               & $ 2.31 \times {10}^{-2} $
               & $ 6.86 \times {10}^{-2} $
                                                                      \\ \hline
               &                              &
               &                                                                                 \\
$ 12$   & $ 1.14 \times {10}^{-2} $
               & $ 1.55 \times {10}^{-2} $
               & $ 2.4 \times {10}^{-3} $
                                                                                     \\ \hline
               &                              &
               &                                                                                 \\
$24$  & $ 1.4 \times {10}^{-3} $
               & $ 4.8 \times {10}^{-3} $
               & $ 6.40 \times {10}^{-4} $
                                                                  \\ \hline
                                                                  &                              &
               &                                                                                 \\
$48$  & $ 2.18 \times {10}^{-4} $
               & $ 1.3 \times {10}^{-3} $
               & $ 2.87 \times {10}^{-4} $
                                                                  \\ \hline
 \end{tabular}
\end{center}
\end{table}
The observed maximum errors
using improved
improved end conditions are summarized in Table 4.
\par \noindent
\begin{table}[htp]
\caption{ Maximum absolute errors for problem $(\ref{e:4.2})$ in
 $ y_{i}.$ }
 \begin{center}
\begin{tabular}{|c|c|} \hline
               &                                           \\
$n$            & $|y(t_{i})-y_{i} |$
                                                            \\ \hline
               &                                            \\
$6 $        & $ 2.53 \times {10}^{-5} $
                                                                                     \\ \hline
               &                                                                               \\
$ 12$  & $ 1.53 \times {10}^{-7} $
                                                                                                    \\ \hline
               &                                                                                \\
$24$   & $ 1.06 \times {10}^{-9}$
                                                                  \\ \hline
                                                                  &                                                                                \\
$48$   & $ 1.09 \times {10}^{-10}$
                                                                  \\ \hline
 \end{tabular}
\end{center}
\end{table}
\par \noindent

\section{N=3: Optimized numerical solution of the 6th order initial value problem}

\subsection{Non-polynomial spline solution for N=3}

A similar optimized numerical method can be developed for the case $N=3$. In this case the
IVP reads

\begin{equation}\label{ee:1.1}
\left. \begin{array}{lcl}
    y^{(6)}(t)+f(t)y(t)&=&g(t),  \ \ \ \ t \ \in \  [a, \ b], \\
    y(a)  \ \ \ = u_{0},\ & &  y^{(1)}(a)   =  u_{1}, \\
    y^{(2)}(a) = u_{2},\ & & y^{(3)}(a)  =   u_{3}, \\
    y^{(4)}(a) = u_{4},\ & & y^{(5)}(a)  =   u_{5},
    \end{array} \right \}
\end{equation}
  Again
  the interval $[a, \ b]$ is divided
  into $n$ equal subintervals, using the grid points $t_{i}  =  a +  ih $
   ($i=0,1,\ldots ,
  n$), where $h \ = \ (b-a)/n.$ \\
    Again we consider the restrictions $S_{i}$ of the solution to each subinterval $[t_{i}, \
  t_{i+1}], \ i=0,1,\ldots , n-1$,
  \begin{eqnarray}\label{ee:1.2}
  S_{i}(t)& = &a_{i}\cos \omega(t-t_{i})+b_{i}\sin
  \omega(t-t_{i})+c_{i}(t-t_{i})^{5}+d_{i}(t-t_{i})^{4} \nonumber \\
  & &
  +e_{i}(t-t_{i})^{3}+q_{i}(t-t_{i})^{2}+r_{i}(t-t_{i})+ v_{i}.
  \end{eqnarray}
     and define
  \begin{equation}\label{ee:1.3}
 \left. \begin{array}{rcl}
   y_{i} \ = \ S_{i}(t_{i}) & & \ \ \ \ M_{i} \ = \ S_{i}^{(2)}(t_{i}),
    \\
 ~~ N_{i} \ = \ S_{i}^{(4)}(t_{i}),  & & \ \ \ \ L_{i} = S_{i}^{(6)}(t_{i})
  ,
  \end{array} \right \} \ i=0,1,\ldots,n.
  \end{equation}
   We denote by $y(t)$ the exact solution of the IVP (\ref{ee:1.1}) and
$y_{i}$ is the approximation to $y(t_{i})$, obtained by the spline
$S(t_{i}).$
  From continuity of the first, third and fifth derivatives at the border points,
   \textit{i.e.}
  $S^{(\mu)}_{i-1}(t_{i})=S^{(\mu)}_{i}(t_{i})$, $\mu=1,3$
  and $5$, one gets

\begin{equation}
   \begin{array} {ll} \label{ee:1.4}
&
h^{6}L_{i-6}\left(\frac{\theta-\sin\theta}{\theta^{6}\sin\theta}-\frac{1}{6\theta^{3}
\sin\theta}+\frac{1}{12\theta\sin\theta}\right)
 \\
 &
+h^{6}L_{i-5}\left(\frac{6}{\theta^{6}}-\frac{2(\cos\theta+2)}{\theta^{5}\sin\theta}+\frac{\cos\theta-1}{3\theta^{3}
\sin\theta} -\frac{\cos\theta-13}{60\theta\sin\theta}\right)
 \\
 &
+h^{6}L_{i-4}
\left(\frac{(8\cos\theta+7)}{\theta^{5}\sin\theta}-\frac{15}{\theta^{6}}+\frac{4\cos\theta+5}{6\theta^{3}
\sin\theta}-\frac{(52\cos\theta-67)}{120\theta\sin\theta}\right)
 \\
&
+h^{6}L_{i-3}\left(\frac{20}{\theta^{6}}-\frac{2(6\cos\theta+4)}{\theta^{5}\sin\theta}-\frac{2(3\cos\theta+1)}{3\theta^{3}
\sin\theta}-\frac{33\cos\theta-13}{30\theta\sin\theta}\right)  \\
&
+h^{6}L_{i-2}\left(\frac{(8\cos\theta+7)}{\theta^{5}\sin\theta}-\frac{15}{\theta^{6}}+\frac{4\cos\theta+5}{6\theta^{3}
\sin\theta}-\frac{(52\cos\theta-67)}{120\theta\sin\theta}\right)  \\
&
+h^{6}L_{i-1}\left(\frac{6}{\theta^{6}}-\frac{2(\cos\theta+2)}{\theta^{5}\sin\theta}+\frac{\cos\theta-1}{3\theta^{3}
\sin\theta} -\frac{\cos\theta-13}{60\theta\sin\theta}\right)  \\
& + h^{6}L_{i}\left(\frac{\theta-\sin\theta}{\theta^{6}\sin\theta}
-\frac{1}{6\theta^{3}
\sin\theta}+\frac{1}{12\theta\sin\theta}\right)
 \\
= &   \left[y_{i-6}  -  6 y_{i-5} +  15 y_{i-4}- 20 y_{i-3}
  +  15y_{i-2}  - 6  y_{i-1} +y_{i}  \right];  \\
  & i=6,7,\ldots,n,
\end{array}
\end{equation}
which can further be written as
\begin{eqnarray}\label{ee:1.5}
& \left(\alpha h^{6}L_{i-6}+\beta h^{6}L_{i-5}+\gamma
h^{6}L_{i-4}+\delta h^{6}L_{i-3}+\gamma h^{6}L_{i-2}+\beta
h^{6}L_{i-1}+\alpha h^{6}L_{i} \right) \nonumber \\
= & \left[y_{i-6}  -  6 y_{i-5} +  15 y_{i-4}- 20 y_{i-3}
  +  15y_{i-2}  - 6  y_{i-1} +y_{i}  \right]; \nonumber \\
  & i=6,7,\ldots,n,
\end{eqnarray}
where
\begin{eqnarray}
\alpha & = & \left(\frac{\theta-\sin\theta}{\theta^{6}\sin\theta}
-\frac{1}{6\theta^{3}
\sin\theta}+\frac{1}{12\theta\sin\theta}\right), \nonumber
\\
\beta & = &
\left(\frac{6}{\theta^{6}}-\frac{2(\cos\theta+2)}{\theta^{5}\sin\theta}+\frac{\cos\theta-1}{3\theta^{3}
\sin\theta} -\frac{\cos\theta-13}{60\theta\sin\theta}\right), \nonumber \\
\gamma & = &
\left(\frac{(8\cos\theta+7)}{\theta^{5}\sin\theta}-\frac{15}{\theta^{6}}+\frac{4\cos\theta+5}{6\theta^{3}
\sin\theta}-\frac{(52\cos\theta-67)}{120\theta\sin\theta}\right)
\nonumber
\end{eqnarray}
and
\begin{eqnarray}
 \delta & = &
\left(\frac{20}{\theta^{6}}-\frac{2(6\cos\theta+4)}{\theta^{5}\sin\theta}-\frac{2(3\cos\theta+1)}{3\theta^{3}
\sin\theta}-\frac{33\cos\theta-13}{30\theta\sin\theta}\right).
\nonumber
\end{eqnarray}
Here $\theta=\omega h.$ The relation $(\ref{ee:1.5})$ forms a system
of $n-5$ linear
equations in the $n$ unknowns \\
$(y_{i}, \ i=1,2,...,n) $, while $L_{i}$ is taken from IVP
$(\ref{e:1.1})$ to be equal to $-f_{i}y_{i} + g_{i}$, \\ $i=0,1,\ldots,n.$ \\
Following~\cite{Sg06}, five equations (end conditions) are
determined to find the complete solution of $y_{i}s$ appearing in
eq. $(\ref{ee:1.5})$, as given below:
\begin{eqnarray}\label{ee:3.1}
 L_{0}  +  L_{4} & = & \frac{1}{h^{6}} \ \left[  \ \frac{2905}{12} y_{0}
   -  336 y_{1} +  126 y_{2} -   \frac{112}{3}y_{3}
   + \frac{21}{4}y_{4}  +  175h y_{0}^{(1)}
     \right. \nonumber \\
    & & \ \ \ \ \ \ \ \ \ \left. +42h^{2}y_{0}^{(2)}- \frac{4}{5} \ h^{6}y_{0}^{(6)}
     \right ], \\
     \label{eee:3.2}
     L_{1}  +  L_{5} & = & \frac{1}{h^{6}} \ \left[  \ \frac{797790}{21983} y_{1}
   -  \frac{1660890}{21983} y_{2} +  \frac{1299060}{21983} y_{3} -   \frac{523110}{21983} y_{4}
          \right. \nonumber \\
    & & \ \ \ \ \ \ \ \ \ \left. + \frac{87150}{21983} y_{5}+  \frac{283500}{21983}h y_{0}^{(1)}
     +   \frac{172620}{21983} \ h^{2} y_{0}^{(2)}
-\frac{40167}{21983}h^{6} y_{1}^{(6)}
     \right ], \\
   \label{ee:3.3}
     L_{2}  +  L_{6} & = & \frac{1}{h^{6}} \ \left[  \ \frac{605725}{22267} y_{2}
   -  \frac{108239440}{1803627} y_{3} +  \frac{1103910}{22267} y_{4} -\frac{446800}{22267}y_{5}
    + \frac{5949805}{1803627}y_{6}
          \right. \nonumber \\
    & & \ \ \ \ \ \ \ \ \ \left. + \frac{675200}{85887}hy_{0}^{(1)} +  \frac{700180}{66801} \ h^{2}y_{0}^{(2)}
    +  \frac{851440}{200403} \ h^{3}y_{0}^{(3)}  \right ]\\
   \label{eee:3.4}
     L_{3}  +  L_{7} & = & \frac{1}{h^{6}} \ \left[  - \frac{670672000}{42346017} y_{3}
   +  \frac{44149995}{1568371} y_{4} -  \frac{23862240}{1568371} y_{5} + \frac{122902615}{42346017}y_{6}
              \right. \nonumber \\
    & & \ \ \ \ \ \ \ \ \ \left. - \frac{12961750}{2016477}hy_{0}^{(1)} -  \frac{25078370}{1568371} \ h^{2}y_{0}^{(2)}
    -  \frac{77684300}{4705113} \ h^{3}y_{0}^{(3)} \right. \nonumber \\
    & & \ \ \ \ \ \ \ \ \ \left. - \frac{11492010}{1568371} \ h^{4}y_{0}^{(4)} \right] \nonumber \\
   \end{eqnarray}
and
\begin{eqnarray}\label{ee:3.5}
  L_{4}  +   L_{8}  & = & \frac{1}{h^{6}} \ \left[ \ \frac{49567095}{12837314}y_{4}
   -  \frac{34289280}{6418657}y_{5} +  \frac{19011465}{12837314} y_{6}
        +  \frac{2182545}{916951}hy_{0}^{(1)}  \right. \nonumber \\
   & & \ \ \ \ \ \ \ \ \  \left.
   +   \frac{59244435}{6418657} \ h^{2}y_{0}^{(2)} + \frac{107795790}{6418657}h^{3}y_{0}^{(3)}
   + \frac{115282605}{6418657}h^{4}y_{0}^{(4)}\right. \nonumber \\
   & & \ \ \ \ \ \ \ \ \  \left. + \frac{65492262}{6418657}h^{5}y_{0}^{(5)}     \right]\ .
    \end{eqnarray}
    The local truncation errors associated with the linear equations
    $(\ref{ee:3.1})-(\ref{ee:3.5})$ and $(\ref{ee:1.5})$ are calculated, as
\begin{equation}\label{ee:3.12}
\tilde{t}_{i} = \left\{ \begin{array}{ll}
   - 4.75 h^{8}y^{(8)}(t_{1})+O(h^{9}), &  i=1,  \\
   - 5.0467 h^{8}y^{(8)}(t_{2})+O(h^{7}), &  i=2,  \\
    - 5.9909 h^{8}y^{(8)}(t_{3})+O(h^{9}), &  i=3,  \\
    - 12.3201 h^{8}y^{(8)}(t_{4})+O(h^{9}), &  i=4,  \\
    - 23.7869 h^{8}y^{(8)}(t_{5})+O(h^{9}), &  i=5,  \\
   (-1+2 \alpha +2 \beta+2\gamma+\delta) h^{6}y^{(6)}(t_{i}) \\
      +\frac{1}{4}(-1+36 \alpha +16 \beta+4\gamma) h^{8}y^{(8)}(t_{i})+O(h^{9}), &  i=6,7,\ldots,n.  \\
        \end{array} \right.
\end{equation}
To make the truncation errors of the system (\ref{ee:1.5}) of order
$h^{8}$, $\alpha$, $\beta$, $\gamma$ and $\delta$ are taken such
that $\alpha+\beta+\gamma+\frac{\delta}{2}=\frac{1}{2}$ and then
\begin{equation}\label{ee:3.13}
\|\tilde{T}\|=ch^{8}R_{3}=O(h^{8}), \ \ \ \ \ \ \ \ R_{3}= \max_{t
\in [a, \ b]} |y^{(8)}(t)|,
\end{equation}\
 where $c$ is a constant. \\

The solution obtained using the system of
linear equations $(\ref{ee:3.1})-(\ref{ee:3.5})$ and $(\ref{ee:1.5})$
in general is second order convergent. Again, however, the order of accuracy of the method can
be improved significantly to $h^{8}$.
 The local truncation error of the system (\ref{ee:1.5}) can be expressed in the following form
\begin{equation}\label{ee:4321}
\tilde{t}_{i} = \left\{ \begin{array}{ll}
      (-1+2 \alpha +2 \beta+2\gamma+\delta) h^{6}y^{(6)}(t_{i})
      +\frac{1}{4}(-1+36 \alpha +16 \beta+4\gamma) h^{8}y^{(8)}(t_{i}) & \\
      +\frac{1}{240}(-7+1620 \alpha +320 \beta+20\gamma) h^{10}y^{(10)}(t_{i})& \\
      + \frac{1}{7560}(-16+15309 \alpha +1344 \beta+21\gamma) h^{12}y^{(12)}(t_{i}) & \\
      +\frac{1}{120960}(-13+39366 \alpha +1536 \beta+6\gamma) h^{14}y^{(14)}(t_{i})& \\
      +\frac{1}{159667200}(-651+5196312 \alpha +90112 \beta+88\gamma) h^{16}y^{(16)}(t_{i}) & \\
      +O(h^{18}), &  \\
      i=6,7,\ldots,n,  & \\
 \end{array} \right.
\end{equation}
Thus, the order of the truncation error $\tilde{t}_{i}$ can be
improved to be of order $h^{14}$ and correspondingly the order of
method can be improved up to $h^{8}$ , if $\alpha=\frac{1}{30240}, \
\beta=41/5040, \ \gamma=2189/10080, \ \delta=\frac{4153}{7560}.$ 
For other choices of the parameters (not listed here), one can make the
method to be of order $h^2,h^4,h^6$, respectively. 
Results corresponding to the order $h^{2}$, $h^{4}$, $h^{6}$ and
$h^{8}$ are described in the following section.

\subsection{Test of the method with analytically solvable examples}
   \ \ \\
{\large {\bf Example 3 }} \\
 \ \ \\
 Consider the IVP
 \begin{equation}
\left.
\begin{array}{ll}\label{ee:4.1}
    y^{(6)}(t) \ -  y(t) \ =  \ -6 e^{t}, \ \ \ \ \ 0\leq t \leq 1 \\
 y(0) \  \ \ \  = 1,  \ \ \ \ \ \ \ \ \ \ \ \ \ \ \ \  y^{(1)}(0) \ = 0, \\
   y^{(2)}(0) \ =  -1, \ \ \ \ \ \ \ \ \ \ \ \ \ y^{(3)}(0) \ =  -2, \\
   y^{(4)}(0) \ =  -3, \ \ \ \ \ \ \ \ \ \ \ \ \ y^{(5)}(0) \ =  -4.
\end{array}
\right\}
\end{equation}
The analytic solution of this IVP $(\ref{ee:4.1})$ is
$$ y(t) \ = \ (1-t) \  e^{t} \ .$$
The observed maximum errors are summarized in Table 5. It is
confirmed from Table 5 that if $h$ is reduced by factor $1/2$, then
$\| E \|$ is reduced by a factor $1/4$, which indicates that the
method gives second-order results.
\par \noindent
\begin{table}[htp]
\caption{ Maximum absolute errors for problem $(\ref{ee:4.1})$ in
 $ y_{i}.$ }
 \begin{center}
\begin{tabular}{|c|c|c|c|} \hline
               &                              &
               &                                           \\
$n$              &  $\alpha=1/120, \ \beta=15/120 $
                &
                  $\alpha=\frac{1}{720}, \ \beta=\frac{1}{36} $ &
                  $\alpha=\frac{1}{5040}, \ \beta=\frac{6}{504} $
                                      \\
               &               &
               &
                                                       \\
               &  $\gamma=1/4,~~~~\delta=28/120 $             &
              $\gamma=\frac{219}{720}, \ \delta=\frac{240}{720} $ &
              $\gamma=\frac{1250}{5040}, \ \delta=\frac{2418}{5040} $
                                                      \\
               &                &
               &
                                                            \\ \hline
               &                              &
               &                                            \\
$8 $        & $ 7.98 \times {10}^{-4} $
               & $ 9.13 \times {10}^{-4} $
               & $ 9.51 \times {10}^{-4} $
                                                                      \\ \hline
               &                              &
               &                                                                                 \\
$ 16$   & $ 7.50 \times {10}^{-5} $
               & $ 9.64 \times {10}^{-5} $
               & $ 1.03 \times {10}^{-4} $
                                                                                     \\ \hline
               &                              &
               &                                                                                 \\
$32$  & $ 5.45 \times {10}^{-6} $
               & $ 1.02 \times {10}^{-5} $
               & $ 1.18 \times {10}^{-5} $
                                                                  \\ \hline
                                                                  &                              &
               &                                                                                 \\
$64$  & $ 1.28 \times {10}^{-7} $
               & $ 9.42 \times {10}^{-7} $
               & $ 1.37 \times {10}^{-6} $
                                                                  \\ \hline
 \end{tabular}
\end{center}
\end{table}
The observed maximum errors (in absolute values) associated with $
y_{i}$, for the problem $(\ref{ee:4.1})$, corresponding to different
orders of method are summarized in Table 6.
\par \noindent
\begin{table}[htp]
\caption{ Maximum absolute errors for problem $(\ref{ee:4.1})$ in
 $ y_{i}.$ }
 \begin{center}
\begin{tabular}{|c|c|c|c|} \hline
               &       & &                                    \\
$n$            & $O(h^{4})$ & $O(h^{6})$ & $O(h^{8})$
                                                            \\ \hline
               &         & &                                    \\
$8 $      & $ 4.04 \times {10}^{-5}$ & $ 2.07 \times {10}^{-1} $  &
$ 2.13 \times {10}^{-1} $
                                                                                     \\ \hline
            & &   &                                                                               \\
$ 16$ &$ 1.10 \times {10}^{-6} $ & $ 8.99 \times {10}^{-9} $ & $
4.80 \times {10}^{-7} $
                                                                                                    \\ \hline
 \end{tabular}
\end{center}
\end{table}
\par \noindent
{\large {\bf Example 4 }} \\
 \ \ \\
 The second example is
 \begin{equation}
\left.
\begin{array}{ll}\label{ee:4.2}
    y^{(6)}(t) \ +  y(t) \ =  \ 6(2t \cos(t)+5 \sin(t)), \ \ \ \ \ -1 \leq t \leq 1 \\
 y(-1) \  \ \ \  = 0,  \ \ \ \ \ \ \ \ \ \ \ \ \ \ \ \  y^{(1)}(-1) \ = 2~sin(1), \\
   y^{(2)}(-1) \ =  -4cos(1)-2sin(1), \ \ \ \ \ \ \ \ \ \ \ \ \ y^{(3)}(-1) \ =  6cos(1)-6sin(1), \\
   y^{(4)}(-1) \ =  8cos(1)+12sin(1), \ \ \ \ \ \ \ \ \ \ \ \ \ y^{(5)}(-1) \ =  -20cos(1)+10sin(1).
\end{array}
\right\}
\end{equation}
with analytic solution
$$ y(t) \ = \ (t^{2}-1) \  \sin(t) \ .$$
The observed maximum errors (in absolute values) associated with $
y_{i}$, for the system $(\ref{ee:4.2})$, corresponding to the
different values of $\alpha, \ \beta$ $\gamma$ and $\delta$ are
summarized in Table 7. Again it is confirmed from Table 7 that if $h$ is
reduced by factor $1/2$, then $\| E \|$ is reduced by a factor
$1/4$, which indicates that the method gives second-order results.
\par \noindent
\begin{table}[htp]
\caption{ Maximum absolute errors for problem $(\ref{ee:4.2})$ in
 $ y_{i}.$ }
 \begin{center}
\begin{tabular}{|c|c|c|c|} \hline
               &                              &
               &                                           \\
$n$              &  $\alpha=1/120, \ \beta=15/120 $
                &
                  $\alpha=\frac{1}{720}, \ \beta=\frac{1}{36} $ &
                  $\alpha=\frac{1}{5040}, \ \beta=\frac{6}{504} $
                                      \\
               &               &
               &
                                                       \\
               &  $\gamma=1/4,~~~~\delta=28/120 $             &
              $\gamma=\frac{219}{720}, \ \delta=\frac{240}{720} $ &
              $\gamma=\frac{1250}{5040}, \ \delta=\frac{2418}{5040} $
                                                      \\
               &                &
               &
                                                            \\ \hline
               &                              &
               &                                            \\
$16 $        & $ 7.35 \times {10}^{-2} $
               & $ 9.64 \times {10}^{-2} $
               & $ 1.03 \times {10}^{-1} $
                                                                      \\ \hline
               &                              &
               &                                                                                 \\
$ 32$   & $ 1.01 \times {10}^{-2} $
               & $ 1.62 \times {10}^{-2} $
               & $ 1.82 \times {10}^{-2} $
                                                                                     \\ \hline
               &                              &
               &                                                                                 \\
$64$  & $ 4.51 \times {10}^{-4} $
               & $ 2.0 \times {10}^{-3} $
               & $ 2.5 \times {10}^{-3} $
                                                                  \\ \hline
                                                                  &                              &
               &                                                                                 \\
$128$  & $ 1.98 \times {10}^{-4} $
               & $ 1.79 \times {10}^{-4} $
               & $ 3.05 \times {10}^{-4} $
                                                                  \\ \hline
 \end{tabular}
\end{center}
\end{table}
The observed maximum errors
corresponding to the
different orders of the method are summarized in Table 8.
\par \noindent
\begin{table}[htp]
\caption{ Maximum absolute errors for problem $(\ref{ee:4.2})$ in
 $ y_{i}.$ }
 \begin{center}
\begin{tabular}{|c|c|c|c|} \hline
               &       & &                                    \\
$n$            & $O(h^{4})$ & $O(h^{6})$ & $O(h^{8})$
                                                            \\ \hline
&         & &                                    \\
$8 $      & $ 2.31 \times {10}^{-2}$ & $ 2.87 \times {10}^{-1} $  &
$ 2.98 \times {10}^{-1} $
                                                                                     \\ \hline
               &         & &                                    \\
$16 $      & $ 8.6 \times {10}^{-3}$ & $ 7.98 \times {10}^{-5} $  &
$ 9.93 \times {10}^{-8} $
                                                                                     \\ \hline
 \end{tabular}
\end{center}
\end{table}
\par \noindent

\par \noindent
\section{Conclusion and Outlook}
In this paper we started from a set of $N$ nonlocal coupled harmonic oscillators,
each driven by a driving force. We showed that this leads to an $2N$-th order
initial value problem (IVP) in a single variable. In a sense this gives `physical meaning'
to high-order IVP in one variable, which so far have mainly been looked at without any
physical interpretation. Engineering applications include strongly coupeld 
oscillator problems where the state of a local oscillator 
is strongly influenced by the
position of the nearest neighbor.

By implementing improved end conditions, a very precise numerical
method could be developed to solve this system numerically. In fact,
we believe it is one of the most precise methods known in the field.
Apparently our results are relevant to find very precise numerical
solution schemes for higher-dimensional differential equations. We
showed that a transformation of an $N$-dimensional 2nd order
differential equation to a 1-dimensional differential equation of
order $2N$ can be highly advantageous from a numerical point of
view. One can implement improved end conditions that allow for a
significant reduction of the error. After the $2N$-th order IVP has
been solved very precisely, the solution can be translated back into
the original physical setting of $N$ nonlocal oscillators.

While we have explicitly worked out the cases $N=2$ and $N=3$, in principle our method can be
extended to higher values of $N$, though the complexity of the formulas used
to minimize the truncation error increases rapidly.

\end{document}